\documentclass[12pt]{article}
\usepackage{amsmath}
\usepackage{amssymb}
\usepackage{amsthm}
\usepackage{extarrows}

\newtheorem{proposition}{Proposition}

\newcommand{\Rm}{\mathbb{R}^m}

\newcommand{\Clm}{\mathcal{C}l_m}

\newcommand{\DD}{\mathcal{D}_{1,k}}
\newcommand{\Do}{\mathcal{D}_{1,2n-1}}
\newcommand{\De}{\mathcal{D}_{1,2n}}

\newcommand{\Sm}{\mathbb{S}^{m-1}}

\newcommand{\SO}{\mathbb{S}^1}
\newcommand{\DDCL}{\mathcal{D}_{1,k}^{C_l}}
\newcommand{\Z}{\mathbb{Z}}
\newcommand{\m}{\textbf{m}}

\begin{document}

\title{Higher order fermionic and bosonic operators on cylinders and Hopf manifolds}

\author{Chao Ding$^1$\thanks{Electronic address:  {\tt dchao@uark.edu}.}, Raymond Walter$^{1,2}$\thanks{Electronic address:  {\tt rwalter@email.uark.edu}; R.W. acknowledges this material is based upon work supported by the National Science Foundation Graduate Research Fellowship Program under Grant No. DGE-0957325 and the University of Arkansas Graduate School Distinguished Doctoral Fellowship in Mathematics and Physics.}, and John Ryan$^1$\thanks{Electronic address: {\tt jryan@uark.edu.}} \\
\emph{\small $^1$Department of Mathematics, University of Arkansas, Fayetteville, AR 72701, USA} \\ 
\emph{\small $^2$Department of Physics, University of Arkansas, Fayetteville, AR 72701, USA}}

\date{}

\maketitle

\begin{abstract}
Higher order higher spin operators are generalizations of $kth$-powers of the Dirac operator. In this paper, we study higher order higher spin operators defined on some conformally flat manifolds, namely cylinders and Hopf manifolds. We will also construct the kernels of these operators on these manifolds.
\end{abstract}

\section{Introduction}
The Rarita-Schwinger operators, as the generalizations of Dirac operators in the higher spin theory, have been studied in many papers in different spaces: for instance, in Euclidean space \cite{B,D,Li2,Ra}, on spheres and real projective spaces \cite{LRV}, and on cylinders and some other conformally flat manifolds \cite{Li1,Li}. The Maxwell operator and the higher spin Laplace operators are the generalizations of the Laplace operator in the higher spin theory in Euclidean space \cite{B1,E,Ro}. Higher order higher spin operators are first studied on the sphere in \cite{Smid} . Then in \cite{Ding}, higher order higher spin operators are introduced for spin 1 and spin $\frac{3}{2}$ in Euclidean space. However, such operators have not been studied on other conformally flat manifolds, such as cylinders and Hopf manifolds.\\
\par
Conformally flat manifolds are manifolds with atlases whose transition functions are M\"{o}bius transformations. They can be constructed by factoring out a subdomain $U$ of either the sphere $\mathbb{S}^m$ or $\Rm$ by a Kleinian group $\Gamma$ of the M\"{o}bius group, where $\Gamma$ acts strongly discontinuously on $U$ and $\Gamma$ is not cyclic. This gives rise to a conformally flat manifold $U\backslash \Gamma$. Cylinders are examples of conformally flat manifolds of type $\Rm\backslash \Z^l$ where $\Z^l$ is an integer lattice and $1\leq l\leq m$.\\
\par
In this paper, we define $k$th order higher spin operators on cylinders and Hopf manifolds, where $k$ is a positive integer. Their fundamental solutions are also constructed by applying translation groups or dilation groups.\\
\par
R. S. Krausshar and J. Ryan \cite{KR1,KR} introduced Clifford analysis on cylinders and tori, making use of the fact that the universal covering space of all of these manifolds is $\Rm$. Then, provided the functions and kernels are $l$-periodic for some $l\in\{1,2,\cdots,m\}$, they used the projection map to obtain the equivalent function or kernel on those manifolds. J. Li, J. Ryan and J. Vanegas \cite{Li1,Li} used similar ideas to construct Rarita-Schwinger operators on cylinders and some other conformally flat manifolds. We follow them by taking the kernel of the higher order higher spin operators in $\Rm$ and using the translation group to construct new kernels that are $l$-fold periodic, which are then projected onto cylinders. In order to prove that the new kernels are well defined, we adapt the Eisenstein series argument developed by Krausshar \cite{K}.\\
\par
\section{Preliminaries}
\subsection{Clifford algebra}\hspace*{\fill}
A real Clifford algebra $Cl_{m}$ can be generated from $\mathbb{R}^m$ by considering the
relationship $$\underline{x}^{2}=-\|\underline{x}\|^{2}$$ for each
$\underline{x}\in \mathbb{R}^m$.  We have $\mathbb{R}^m\subseteq Cl_{m}$. If $\{e_1,\ldots, e_m\}$ is an orthonormal basis for $\mathbb{R}^m$, then $\underline{x}^{2}=-\|\underline{x}\|^{2}$ tells us that $$e_i e_j + e_j e_i= -2\delta_{ij},$$ where $\delta_{ij}$ is the Kronecker delta function. An arbitrary element of the basis of the Clifford algebra can be written as $e_A=e_{j_1}\cdots e_{j_r},$ where $A=\{j_1, \cdots, j_r\}\subset \{1, 2, \cdots, m\}$ and $1\leq j_1< j_2 < \cdots < j_r \leq m.$
Hence for any element $a\in Cl_{m}$, we have $a=\sum_Aa_Ae_A,$ where $a_A\in \mathbb{R}$.\\
\par
The Dirac operator in $\mathbb{R}^m$ is defined to be $$D_x:=\sum_{i=1}^{m}e_i\partial_{x_i}.$$  Note $D_x^2=-\Delta_x$, where $\Delta_x$ is the Laplacian in $\mathbb{R}^m$.  A $\Clm$-valued function $f(x)$ defined on a domain $U$ in $\Rm$ is left monogenic if $D_xf(x)=0.$ Since multiplication of Clifford numbers is not commutative in general, there is a similar definition for right monogenic functions.

\section{\textbf{The higher order higher spin operator on cylinders }}
For an integer $l$, $1\leq l\leq m$, we define the $l$-cylinder $C_l$ to be the $m$-dimensional manifold $\Rm/\Z^l$, where $\Z^l$ denote the $l$-dimensional lattice defined by $\Z^l:=\Z e_1+\cdots+\Z e_l$. We denote its members $m_1e_1+\cdots+m_le_l$ for each $m_1,\cdots,m_l\in\Z$ by a bold letter \m. When $l=m$, $C_l$ is the $m$-torus, $T_m$. For each $l$ the space $\Rm$ is the universal covering group of the cylinder $C_l$. Hence there is a projection map $\pi_l:\ \Rm\longrightarrow C_l$.\\
\par
An open subset U of the space $\Rm$ is called $l$-fold periodic if for each $x\in U$ the point $x+\m\in U$. So $U':=\pi_l(U)$ is an open subset of the $l$-cylinder $C_l$.\\
\par
Suppose that $U\in\Rm$ is a $l$-fold periodic open set. Let $f(x,u)$ be a function defined on $U\times\Rm$ with values in $\Clm$. Then we say that $f(x,u)$ is a $l$-fold periodic function if for each $x\in U$ we have that $f(x,u)=f(x+\m,u)$. Moreover, we will assume throughout that $f$ is a monogenic polynomial homogeneous of degree $j$ in $u$.\\
\par
Now, if $f:\ U\times\Rm\longrightarrow\Clm$ is an $l$-fold periodic function, then the projection $\pi_l$ induces a well defined function
$$f':\ U'\times\Rm\longrightarrow\Clm,$$
where $f'(x',u)=f(x,u)$ for each $x'\in U'$ and $x$ an arbitrary representative of $\pi_l^{-1}(x')$. Moreover, any function $f':\ U'\times\Rm\longrightarrow\Clm$ lifts to an $l$-fold periodic function $f:\ U\times\Rm\longrightarrow\Clm$, where $U=\pi_l^{-1}(U')$.\\
\par
The projection map $\pi_l$ induces a projection of the higher order higher spin operator $\DD$ to an operator $\DDCL$ acting on domains on $C_l\times\Rm$ which is defined by replacing $D_x$ and $\Delta_x$ with $D'_x$ and $\Delta_x'$ in $\DD$, respectively, where $D'_x$ is the projection of the Dirac operator $D_x$ and $\Delta_x'$ is the projection of the Laplace operator $\Delta_x$. Essentially this is the method in \cite{Li1}. That is
\begin{eqnarray*}
\De^{C_l}&=&\Delta_x'^n-\frac{4n}{m+2n-2}\langle u,D'_x\rangle \langle D_u,D_x'\rangle\Delta_x'^{n-1};\\
\Do^{C_l}&=&D_x'\Delta_x'^{n-1}-\frac{2}{m+2n-2}u\langle D_u,D_x'\rangle \Delta_x'^{n-1}-\frac{4n-4}{m+2n-2}\langle u,D_x'\rangle \langle D_u,D_x'\rangle \Delta_x'^{n-2}D_x'.
\end{eqnarray*}
We call the operator $\De^{C_l}$ an $l$-cylindrical higher order bosonic operator of spin 1 and $\Do^{C_l}$ an $l$-cylindrical higher order fermionic operator of spin $\frac{3}{2}$.
\subsection*{Fundamental solutions of $\DDCL$}
We follow the techniques in \cite{Li1} and \cite{Ding}, requiring that the order of the operator $k<m$ when the dimension $m$ is even. Let $U$ a domain in $\Rm$. We recall the fundamental solution of the higher order higher spin operators $\DD$ in $\Rm$:
\begin{equation*}
E_{1,k}(x,u,v)=c_1G_k(x)Z_1(\displaystyle\frac{xux}{||x||^2},v),
\end{equation*}
where $c_1$ is a non-zero real constant, $Z_1(u,v)$ is the reproducing kernel of degree-$1$ homogeneous harmonic (respectively monogenic) polynomials if $k$ is even (resp. odd), and
\begin{equation*}
G_k(x):=
\begin{cases}
\displaystyle\frac{1}{||x||^{m-2n}},  &\text{if $k=2n$;}\\
\displaystyle\frac{x}{||x||^{m-2n+2}},   &\text{if $k=2n-1$.}
\end{cases}
\end{equation*}\\
\par
Now we consider sums of the following form:
\begin{eqnarray}\label{1}
cot_{l,1,k}(x,u,v)=\sum_{(m_1,\cdots,m_l)\in\Z^l}E_{1,k}(x+m_1e_1+\cdots+m_le_l,u,v),
\end{eqnarray}
for all $1\leq p\leq m-k-1$.\\
\par
We will show these functions are defined on the $l$-fold periodic domain $\Rm/\Z^l$ for fixed $u$ and $v$ in $\Rm$ and are $\mathcal{C}l_m$-valued. They are also $l$-fold periodic functions.\\
\par
To prove the locally uniform convergence of the series (\ref{1}), use the locally normal convergence of the series
$$\sum_{\textbf{m}\in\Z^l}G_k(x+\m)$$
established by the following proposition \cite{K}.
\begin{proposition}
Let $p\in\mathbb{N}$ with $1\leq p\leq m-k-1$. Let $\Z^p$ be the $p$-dimensional lattice. Then the series
\begin{eqnarray*}
\sum_{\m\in\Z^p}q_0^{k}(x+\m)
\end{eqnarray*}
converges normally in $\Rm\backslash\Z^p$.
\end{proposition}
Here $q_0^{k}$ is in the kernel of the $k$th-power of $D_x$ and $q_0^{k}$ is exactly our $G_k$ above. The proof follows from Proposition 2.2 appearing in \cite{K}.\\
\par
Returning to the series defined by (\ref{1}),
\begin{eqnarray*}
&&cot_{l,1,k}(x,u,v)=\sum_{(m_1,\cdots,m_l)\in\Z^l}E_{1,k}(x+m_1e_1+\cdots+m_le_l,u,v)\\
&=&\sum_{\m\in\Z^l}G_k(x+\m)Z_1(\frac{(x+\m)u(x+\m)}{||x+\m||^2},v),\ 1\leq l\leq n-k-1,
\end{eqnarray*}
we observe that $Z_1(\displaystyle\frac{(x+\m)u(x+\m)}{||x+\m||^2},v)$ is a bounded function on a bounded domain in $\Rm$ because its first variable,
$\displaystyle\frac{(x+\m)u(x+\m)}{||x+\m||^2}$, is a reflection in
the direction $\displaystyle\frac{(x+\m)}{||x+\m||}$ for each $\m$, and hence is a linear transformation which
is a continuous function. Furthermore, the norm of $\displaystyle\frac{(x+\m)u(x+\m)}{||x+\m||^2}$ is equal to the norm of $u$, so the bound of $Z_1$ with respect to the first variable does not depend on $\m$. On the other hand, we would get with respect to the second variable, bounded homogeneous functions of degree $1$.\\
\par
Consequently, applying the former proposition, the series
\begin{eqnarray*}
cot_{l,1,k}(x,u,v)=\sum_{(m_1,\cdots,m_l)\in\Z^l}E_{1,k}(x+m_1e_1+\cdots+m_le_l,u,v),
\end{eqnarray*} 
for $1\leq l\leq m-k-1$, is a uniformly convergent series and represents a kernel for the higher order higher spin operators under translations by $\m\in\Z^l$, with $1\leq l\leq m-k-1$.\\
\par
Using similar argument in \cite{Li1}, we define the $(m-k)$-fold periodic cotangent. In order to do that, we decompose the lattice $\Z^l$ into three parts: the origin $\{0\}$ and a positive and a negative part. The last two parts are qual and disjoint:
\begin{eqnarray*}
\Lambda_l&=&\{m_1e_1:\ m_1\in\mathbb{N}\}\cup\{m_1e_1+m_2e_2:\ m_1,m_2\in\Z,m_2>0\}\\
&&\cup\cdots\cup\{m_1e_1+\cdots+m_le_l:\ m_1,\cdots,m_l\in\Z,m_l>0\}
\end{eqnarray*}
and
$$-\Lambda_l=(\Z^l\backslash\{0\})\backslash\Lambda_l.$$
For $l=m-k$, we define
\begin{eqnarray*}
cot_{m-k,1,k}(x,u,v)=E_{1,k}(x,u,v)+\sum_{\m\in\Lambda_{m-k}}\big[E_{1,k}(x+\m,u,v)+E_{1,k}(x-\m,u,v)\big].
\end{eqnarray*}
\par
The following proposition \cite{K} and a similar argument as in \cite{Li1} shows the above series converges normally.
\begin{proposition}
Let $\Z^{m-k}$ be the $(m-k)$-dimensional lattice. Then the series
$$q_0^k(x)+\sum_{\m\in\Z^{m-k}\backslash\{0\}}\big(q_0^k(x+\m)-q_0^k(\m)\big)$$
\end{proposition}
The proof follows Proposition 2.2 appearing in \cite{K}. Hence, it is a kernel for the higher order higher spin operator under translation by $\m\in\Lambda_{m-k}.$\\
\par
For $x,y\in\Rm$, the function $cot_{l,1,k}(x-y,u,v)$ induces functions on $C_l$:
$$cot'_{l,1,k}(x',y',u,v)=cot_{l,1,k}(x-y,u,v).$$
for each $x',y'\in U'$ and $x, y$ arbitrary representatives of $\pi_l^{-1}(x')$ and $\pi_l^{-1}(y')$. These functions are defined on $(C_l\times C_l)\backslash diagonal(C_l\times C_l)$ for each fixed $u,v\in\Rm$, where
$$diagonal(C_l\times C_l)=\{(x',x'):\ x'\in C_l\}$$
 and they satisfy
 $$\DD^{C_l}cot'_{l,1,k}(x',y',u,v)=0.$$
 \subsection*{Conformally inequivalent spinor bundles on $C_l$}
Previously the spinor bundle over $C_l$ was trivial, $C_l\times\Clm$. However, there are $2^l$ spinor bundles on $C_l$. See more details in \cite{Li}. The following construction is for some of the spinor bundles over $C_l$ and all the others can be constructed similarly.\\
\par
First let $p$ be an integer in the set $\{1,2\cdots,l\}$ and consider the lattice $\Z^p:=\Z e_1+\cdots+\Z e_p$. We also consider the lattice $\Z ^{l-p}:=\Z e_{p+1}+\cdots+\Z e_l$. In this case $\Z^l=\{\bold{m}+\bold{n}:\m\in\Z^p,\ \bold{n}\in\Z^{l-p}\}$. Suppose that $\m=m_1e_1+\cdots+m_pe_p$. Let us make the identification $(x,X)$ with $(x+\bold{m}+\bold{n},(-1)^{m_1+\cdots+m+p}X)$ where $x\in\Rm$ and $X\in\Clm$. This identification gives rise to a spinor bundle $E^p$ over $C_l$.\\
\par
We adapt functions in previous section as follows. For $1\leq l\leq m-k-1$ we define
\begin{eqnarray*}
cot_{l,1,k,p}(x,u,v)=\sum_{m\Z^p,\ n\in \Z^{l-p}}(-1)^{m_1+\cdots+m_p}E_{1,k}(x+\bold{m}+\bold{n},u,v).
\end{eqnarray*}
These are well defined functions on $\Rm\backslash \Z^l$. Therefore, we obtain from these functions the cotangent kernels
\begin{eqnarray*}
cot_{l,1,k,p}(x,y,u,v)=\sum_{m\Z^p,\ n\in \Z^{l-p}}(-1)^{m_1+\cdots+m_p}E_{1,k}(x-y+\bold{m}+\bold{n},u,v).
\end{eqnarray*}
Again applying the projection map $\pi_l$ these kernels give rise to the kernels
$$cot'_{l,1,k,p}(x',y',u,v).$$
In the case $l=m-k$, by considering the series
\begin{eqnarray*}
&&cot_{m-k,1,k}(x,y,u,v)=E_{1,k}(x,u,v)+\sum_{\m\in\Lambda_{m-k}}\big[E_{1,k}(x+\m,u,v)+E_{1,k}(x-\m,u,v)\big],
\end{eqnarray*}
we obtain the kernel
\begin{eqnarray*}
&&cot_{m-k,1,k}(x,y,u,v)\\
&=&E_{1,k}(x-y,u,v)+\sum_{\m\in\Lambda_{m-k}}\big[E_{1,k}(x-y+\m,u,v)+E_{1,k}(x-y-\m,u,v)\big],
\end{eqnarray*}
which in turn using the projection map induces kernels
$$cot'_{m-k,1,k}(x',y',u,v).$$
Defining
\begin{eqnarray*}
&&cot_{m-k,1,k,p}(x,u,v)=E_{1,k}(x+\bold{m}+\bold{n},u,v)\\
&&+\sum_{\substack{\m\in\Z^p,\bold{n}\in\Z^{m-k-p}\\ \bold{m}+\bold{n}\in\Lambda_{m-k}}}(-1)^{m_1+\cdots+m_p}\big[E_{1,k}(x+\bold{m}+\bold{n},u,v)+E_{1,k}(x-\bold{m}-\bold{n},u,v)\big],
\end{eqnarray*}
we obtain the cotangent kernels
\begin{eqnarray*}
&&cot_{m-k,1,k,p}(x,y,u,v)=E_{1,k}(x-y+\bold{m}+\bold{n},u,v)\\
&&+\sum_{\substack{\m\in\Z^p,\bold{n}\in\Z^{m-k-p}\\ \bold{m}+\bold{n}\in\Lambda_{m-k}}}(-1)^{m_1+\cdots+m_p}\big[E_{1,k}(x-y+\bold{m}+\bold{n},u,v)+E_{1,k}(x-y-\bold{m}-\bold{n},u,v)\big],
\end{eqnarray*}
and by $\pi_l$ the kernels
$$cot'_{m-k,1,k,p}(x',y',u,v).$$
\section{\textbf{The higher order higher spin operator on Hopf manifolds }}
In this section, we use similar arguments as in \cite{Li} and \cite{Ding}, also requiring that the order of the operator satisfies $k<m$ when the dimension $m$ is even. It is worth pointing out that the two spin structures introduced there also apply for our fermionic and bosonic operators. Let $U=\Rm$ and $\Gamma=\{t^i:\ i\in\Z\}$, where $t$ is an arbitrary strictly positive real number distinct from $1$. Then by factoring out $U$ by $\Gamma$ we have the conformally flat spin manifold $\SO\times\Sm$ which we denote by $H_m$ and call a Hopf manifold. As $\prod_1(\SO\times\Sm)=\Z$ for $m>2$, it follows that the Hopf manifold $H_m$ has two distinct spin structures. On this space, we can define the higher order higher spin operators and construct their kernels over the two different spinor bundles over $\SO\times\Sm$.\\
\par
In all that follows $\Rm\backslash\{0\}$ will be a universal covering space of the conformally flat manifold $\SO\times\Sm$ (\cite{KR}). So there is a projection map $p:\ \Rm\backslash\{0\}\longrightarrow\SO\times\Sm$. Further for each $x\in\Rm\backslash\{0\}$ we shall denote $p(x)$ by $x'$. Further if $V$ is a subset of $\Rm\backslash\{0\}$ then we denote $p(V)$ by $V'$.\\
\par
One spinor bundle $F_1$ over $\SO\times\Sm$ can be constructed by identifying the pair $(x,X)$ with $(t^ix,t^{\frac{i(m-1)}{2}}X)$ for every $k\in\Z$, where $x\in\Rm\backslash\{0\}$ and $X\in\mathcal{C}l_m$.\\
\par
Consider a domain $V\subset\Rm\backslash\{0\}$ satisfying $t^ix\in V$ for each $i\in\Z$ and $x\in V$. We will call a such domain a $k$-factor dilation domain. Further we define a $i$-factor dilation function as a function $f(x,u):\ V\times\Rm\longrightarrow\Clm$ such that $f$ is a monogenic polynomial homogeneous of degree $l$ in $u$ satisfying $f(x,u)=t^{\frac{i(m-1)}{2}}f(t^ix,u)$ for each $x\in V$ and each integer $i$.\\
\par
The projection map $p$ induces a well defined function
$$f':\ V'\times\Rm\longrightarrow F_1,$$
where $f'(x,u)=t^{\frac{i(m-1)}{2}}f(x,u)$ for each $x'\in V'$ and $x$ an arbitrary representative of $\pi_l^{-1}(x')$.\\
\par
The higher order higher spin operator over $\Rm\backslash\{0\}$ induces a higher order higher spin operator acting on sections of the bundle $F_1$ over $H_m$. We will denote this operator by $\DD^{H_m}$. If $\DD^{H_m}(f')=0$ then $f'$ is called an $F_1$-left higher order higher spin section. Moreover, any $F_1$-left higher order higher spin section $f':\ V'\times\Rm\longrightarrow F_1$ lifts to a $k$-factor dilation function $f:\ V\times\Rm\longrightarrow \Clm$, where $V=p^{-1}(V')$ and $\DD(f)=0$.\\
\par
Now we consider the series
\begin{eqnarray*}
&&E_{1,k,1}^{H_m}(x,y,u,v)=\sum_{i=-\infty}^0G_k(t^ix-t^iy)Z_1(\frac{(t^ix-t^iy)u(t^ix-t^iy)}{||(t^ix-t^iy)||^2},v)\\
&+&t^{2(k-m)}G_k(x)Z_1(\frac{xux}{||x||^2},v)\big[\sum_{i=1}^{\infty}G_k(t^{-i}x^{-1}-t^{-i}y^{-1})\\
&& Z_1(\frac{(t^{-i}x^{-1}-t^{-i}y^{-1})u(t^{-i}x^{-1}-t^{-i}y^{-1})}{||t^{-i}x^{-1}-t^{-i}y^{-1}||^2},v)\big]G_k(y)Z_1(\frac{yuy}{||y||^2},v),
\end{eqnarray*}
where $x,y\in\Rm\backslash\{0\}$ and $y\neq t^ix$ for all $i\in\Z$. From the definition of $Z_1(u,v)$ and the homogeneity of $G_k(x)$, it is easy to obtain that $E_{1,k,1}^{H_m}$ converges normally on any compact subset $K$ not containing the points $y=t^ix$ where $i\in\Z$. Since $E_{1,k,1}^{H_m}(tx,ty,u,v)=E_{1,k,1}^{H_m}(x,y,u,v)$, this kernel is periodic with respect to the Kleinian group $\{t^i:\ i\in\Z\}$. The kernel for the higher order higher spin operators on $(\SO\times\Sm)\times(\SO\times\Sm)\backslash diagonal(\SO\times\Sm)$ is then the projection of $E_{1,k,1}^{H_m}(x,y,u,v)$ on $(\SO\times\Sm)\times(\SO\times\Sm)\backslash diagonal(\SO\times\Sm)$:
$$E_{1,k,1}^{H_m}(x',y',u,v)=E_{1,k,1}^{H_m}(x,y,u,v),$$
for $x,y$ representatives of $\pi_l^{-1}(x')$ and $\pi_l^{-1}(y')$ as earlier.
\par
The second spinor bundle $F_2$ over $\SO\times\Sm$ can be constructed by identifying the pair $(x,\ X)$ with $(t^ix,(-1)^it^{\frac{i(m-1)}{2}}X)$.\\
\par
Now we introduce the normally convergent series:
\begin{eqnarray*}
&&E_{1,k,2}^{H_m}(x,y,u,v)=\sum_{i=-\infty}^0(-1)^iG_k(t^ix-t^iy)Z_1(\frac{(t^ix-t^iy)u(t^ix-t^iy)}{||(t^ix-t^iy)||^2},v)\\
&+&(-1)^it^{2(k-m)}G_k(x)Z_1(\frac{xux}{||x||^2},v)\big[\sum_{i=1}^{\infty}G_k(t^{-i}x^{-1}-t^{-i}y^{-1})\\
&& Z_1(\frac{(t^{-i}x^{-1}-t^{-i}y^{-1})u(t^{-i}x^{-1}-t^{-i}y^{-1})}{||t^{-i}x^{-1}-t^{-i}y^{-1}||^2},v)\big]G_k(y)Z_1(\frac{yuy}{||y||^2},v)
\end{eqnarray*}
where $x,y\in\Rm\backslash\{0\}$ and $y\neq t^ix$ for all $i\in\Z$. This function induces through the projection map $p$ on the variable $x,y\in\Rm\backslash\{0\}$, the higher order higher spin kernel associated with $F_2:$
\begin{eqnarray*}
E_{1,k,2}^{H_m}(x',y',u,v)=E_{1,k,2}^{H_m}(x,y,u,v),
\end{eqnarray*}
for $x$ and $y$ again as earlier.

\end{document}